\documentclass[12pt]{amsart}

\numberwithin{equation}{section}

\newtheorem{theorem}[equation]{Theorem}
\newtheorem{lemma}[equation]{Lemma}
\newtheorem{corollary}[equation]{Corollary}
\newtheorem{proposition}[equation]{Proposition}
\newtheorem{multonelemma}[equation]{Multiplicity One Lemma}
\newtheorem{mainconj}[equation]{The Main Conjecture}

\theoremstyle{definition}
\newtheorem{definition}[equation]{Definition}
\newtheorem{example}[equation]{Example}

\theoremstyle{remark}

\newenvironment{pf}
{\noindent\textbf{Proof: }}
{\hfill\textbf{ Q.E.D.}\medskip}

\newcommand{\bP}{{\mathbb{P}}}
\newcommand{\bF}{{\mathbb{F}}}
\renewcommand{\L}{\mathcal{L}}
\newcommand{\M}{\mathcal{M}}
\newcommand{\Lo}{\mathcal{L}_0}
\newcommand{\Lohat}{\hat{\mathcal{L}}_0}
\newcommand{\LP}{\mathcal{L}_\bP}
\newcommand{\LF}{\mathcal{L}_\bF}
\newcommand{\LPhat}{\hat{\mathcal{L}}_\mathbb{P}}
\newcommand{\LFhat}{\hat{\mathcal{L}}_\mathbb{F}}

\newcommand{\RP}{\mathcal{R}_\mathbb{P}}
\newcommand{\RF}{\mathcal{R}_\mathbb{F}}
\renewcommand{\l}{\ell}
\newcommand{\lo}{\ell_0}
\newcommand{\lohat}{\hat{\ell}_0}
\newcommand{\lP}{\ell_\mathbb{P}}
\newcommand{\lF}{\ell_\mathbb{F}}
\newcommand{\lPhat}{\hat{\ell}_\mathbb{P}}
\newcommand{\lFhat}{\hat{\ell}_\mathbb{F}}
\newcommand{\rP}{r_\mathbb{P}}
\newcommand{\rF}{r_\mathbb{F}}
\renewcommand{\v}{v}

\newcommand{\vP}{v_\mathbb{P}}
\newcommand{\vF}{v_\mathbb{F}}
\newcommand{\vPhat}{\hat{v}_\mathbb{P}}
\newcommand{\vFhat}{\hat{v}_\mathbb{F}}
\newcommand{\e}{e}
\renewcommand{\O}{\mathcal{O}}
\newcommand{\Pic}{\mathrm{Pic}}

\title{Degenerations of Planar Linear Systems}
\author{Ciro Ciliberto and Rick Miranda}
\address{Dept. of Mathematics\\Universit\'a di Roma II\\Via Fontanile di
Carcaricola\\00173 Rome, Italy}
\email{cilibert@axp.mat.utovrm.it}
\address{Dept. of Mathematics\\Colorado State University\\Ft. Collins, CO
80523}
\email{miranda@math.colostate.edu}
\thanks{Research supported in part by the NSA}
\thanks{\today. Set in type by LaTeX2e}

\begin{document}

\begin{abstract}
Fixing $n$ general points $p_i$ in the plane,
what is the dimension of the space of plane curves of degree $d$
having multiplicity $m_i$ at $p_i$ for each $i$?
In this article we propose an approach to attack this problem,
and demonstrate it by successfully computing this dimension
for all $n$ and for $m_i$ constant, at most $3$.
This application,
while previously known (see \cite{hirschowitz1}),
demonstrates the utility of our approach,
which is based on an analysis of the corresponding linear
system on a degeneration of the plane itself,
leading to a simple recursion for these dimensions.
We also obtain results in the ``quasi-homogeneous'' case
when all the multiplicities are equal except one; this is the
natural family to consider in the recursion.
\end{abstract}

\maketitle

\tableofcontents

\section*{Introduction}
Fix the projective plane $\bP^2$
and $n+1$ general points $p_0, p_1, \ldots, p_n$ in it.
Let $H$ denote the line class of the plane.
Consider the linear system consisting of
plane curves of degree $d$ (that is, divisors in $|dH|$)
with multiplicity $m_0$ at $p_0$
and multiplicity $m_i$ at $p_i$ for $i \geq 1$.
If all $m_i$ for $i \geq 1$ are equal, to $m$ say,
we denote this system by $\L = \L(d,m_0,n,m)$
and call the system \emph{quasi-homogeneous}.
Define its {\em virtual dimension}
\[
\v= \v(d,m_0,n,m) = d(d+3)/2 - m_0(m_0+1)/2 - nm(m+1)/2);
\]
the linear system of all plane curves of degree $d$
has dimension $d(d+3)/2$
and a point of multiplicity $k$ imposes $k(k+1)/2$ conditions.
Of course the actual dimension of the linear system
cannot be less than $-1$
(projectively, dimension $-1$ means an empty system);
hence we define the {\em expected dimension} to be
\[
\e = \e(d,m_0,n,m) = \max\{-1, \v(d,m_0,n,m)\}.
\]

As the points $p_i$ vary on $\bP^2$,
the dimension of this linear system is upper semi-continuous;
therefore on a Zariski open set in the parameter space of
$(n+1)$-tuples of points, the dimension achieves its minimum value,
which we (abusing notation slightly) call the {\em dimension} of $\L$,
and denote by $\l = \l(d,m_0,n,m)$.
We always have that
\[
\l(d,m_0,n,m) \geq \e(d,m_0,n,m);
\]
equality implies
(when the numbers are at least $-1$)
that the conditions imposed by the multiple points are independent.

We will say that the generic system $\L$ is {\em non-special}
if equality holds, i.e.,
that either the system is empty
or that the conditions imposed by the multiple points are independent.
If $\l > \e$ then we say the system is {\em special}.

In this article we discuss the speciality
of such linear systems of plane curves,
and classify the special linear systems
with low multiplicity $m$.

There is a long history for this problem;
we will not make an attempt here to review it.
The ``homogeneous'' cases with $m_0=0$,
and $m \leq 2$, are discussed in \cite{arbarello-cornalba},
and in \cite{hirschowitz1}; this last reference
also has results on the homogeneous $m=3$ case.
In \cite{alexander}, \cite{alex-hirsch1},
\cite{alex-hirsch2}, \cite{evain}, \cite{gergimhar}, \cite{gergimpit},
\cite{harbourne1}, \cite{harbourne2}, \cite{harbourne3},
\cite{hirschowitz3}, \cite{mignon}, and \cite{segre}
one may also find related conjectures
and results on the general case.
The reader may consult
\cite{gimigliano} for a survey.

In Section 1 we lay out some basic notation and elementary observations.
In Section 2 we describe in detail our approach,
which is based on a degeneration of the plane
and of the corresponding linear system.
This leads to a recursion for the sought-after dimension,
which relies on a transversality theorem
for the pair of linear systems
to which the recursion is reduced;
this is described in Section 3.
The failure of systems to be non-special is always due
(in our examples) to the presence of multiple $(-1)$-curves
in the base locus;
we formalize this in Section 4
(calling such systems ``$(-1)$-special'')
and give a classification of quasi-homogeneous $(-1)$-curves
in Section 5.
In Section 6 we present a computation of the dimension
of the linear system $\L(d,m_0,n,m)$ with large $m_0$,
which is particularly useful in our inductive approach;
this uses Cremona transformations in an essential way.
Using this, we present a list of the $(-1)$-special systems
in Section 7, for $m \leq 3$.
We then prove in Section 8
that all quasi-homogeneous special systems with $m \leq 3$
are $(-1)$-special.

In the paper \cite{cm2} we turn our attention to the homogeneous case.
We classify all $(-1)$-special homogeneous systems,
and prove every special homogeneous system with $m \leq 12$
is $(-1)$-special. The main technique which is used for this
is the one described in the present paper.

The authors would like to thank L. Caporaso, A. Geramita,
B. Harbourne, J. Harris, and A. Tjurin
for some useful discussions.

\section{Basic Facts}

We keep the notation of the Introduction.
If there is no danger of confusion,
we will omit in what follows the indication of the data $(d,m_0,n,m)$.
Also, if $n=0$, we omit $n$ and $m$ from the notation
and speak of the linear system $\L(d,m_0)$.

We note that the speciality of $\L$
is equivalent to a statement about linear systems on
the blowup $\bF_1$ of $\bP^2$ at the point $p_0$.
If $H$ denotes the class of the pullback of a line
and $E$ denotes the class of the exceptional divisor,
then the linear system on $\bP^2$
transforms to the linear subsystem of $|dH-m_0E|$
consisting of those curves having multiplicity $m$
at each of the $n$ transformed points $p_i, i \geq 1$
(none of which lie on $E$).
If we further blow up the points $p_i$ to $E_i$,
we obtain the rational surface $\bP^\prime$;
the linear system of proper transforms
is the complete linear system $|dH - m_0E-\sum_i m_iE_i|$
and we denote by $\L^\prime = \L^\prime(d,m_0,n,m)$
the corresponding line bundle.
Then the original system $\L$ is non-special if and only if
\[
h^1(\L^\prime) = \max\{0, -1-\v \}.
\]
In particular if the system is non-empty
(which means that $H^0(\L^\prime)$ is non-zero),
it is non-special if and only if the $H^1$ is zero.
More precisely if the virtual dimension $\v \geq -1$
then non-speciality means that the $H^1$ is zero,
or, equivalently, that the conditions imposed by
the multiple base points are linearly independent.

There is one situation in which a completely general statement
can be made, namely the case of simple base points.
Suppose that $\M$ is a linear system on a variety $X$,
and $p_1,\ldots,p_n$ are general points of $X$.
Let $\M(-\sum_i p_i)$ be the linear subsystem
consisting of those divisors in $\M$
which pass through the $n$ given points.

\begin{multonelemma}
\label{multone}
$\dim \M(-\sum_{i=1}^n p_i) = \max\{-1, \dim\M - n\}$.
\end{multonelemma}

\begin{pf}
The proof is obtained by induction on $n$;
one simply chooses the next point not to be a base
point of the previous linear system, if this is non-empty.
\end{pf}

We can speak of the
\emph{self-intersection} ${\L}^2$ and the
\emph{genus} $g_{\L}$
of the curves of the system ${\L}(d,m_0,n,m)$,
which will be the self-intersection and the arithmetic genus
of $\L^\prime$ on the blowup $\bP^\prime$.
We have:

\[
{\L}^2=d^2-m^2_0-nm^2
\;\;\;\mathrm{ and }\;\;\;
2g_{\L}-2=d(d-3)-m_0(m_0-1)-nm(m-1).
\]
Notice the basic identity:

\begin{equation}
\label{RR}
{\v}={\L}^2-g_{\L}+1.
\end{equation}

We can also speak of the \emph{intersection number}
${\L}(d,m_0,n,m)\cdot{\L}(d',m_0',n',m')$, $n'\leq n$,
which is of course given by:
\[
{\L}(d,m_0,n,m)\cdot{\L}(d',m'_0,n',m) := dd' - m_0m_0' - n'mm'
\]

We collect here certain initial observations.

\begin{lemma}\label{basic1}\mbox{}
\begin{itemize}
\item[a.] $\L(d,0,1,m)$ is non-special for all $d$, $m$.
\item[b.] If $\e(d,m_0,n,m)\geq -1$ and $\L(d,m_0,n,m)$ is non-special
then $\L(d',m_0',n',m')$ is non-special
whenever $d' \geq d$, $m_0' \leq m_0$, $n' \leq n$, and $m' \leq m$.
\item[c.] If $\e(d,m_0,n,m) = -1$ and $\L(d,m_0,n,m)$ is non-special
(i.e., the system $\L$ is empty)
then $\L(d',m_0',n',m')$ is non-special (and therefore empty)
whenever $d' \leq d$, $m_0' \geq m_0$, $n' \geq n$, and $m' \geq m$.
\item[d.] If $d \geq 2$ then $\L(d,0,2,d)$ and $\L(2d,0,5,d)$
are special.
\end{itemize}
\end{lemma}

\begin{pf}
Statement (a) is elementary.

For statement (b), we note that if $d'=d$,
then the inequalities on $m_0'$, $n'$, and $m'$
imply that the conditions imposed on the curves of degree $d$
for the system $\L(d,m_0',n',m')$
are a subset of the conditions
for the system $\L(d,m_0,n,m)$;
since those are independent by the assumptions
that $\L$ is non-special and non-empty,
so are the conditions in the subset.
and therefore $\L(d,m_0',n',m')$ is non-special.
Therefore to prove (b) we may assume that
$m_0'=m_0$, $n'=n$, and $m'=m$, and using induction that $d'=d+1$.
If we now pass to the blow-up ${\bP}'$ of ${\bP}^2$,
we have the following exact sequence:
\[
0 \to {\L}^\prime(d,m_0,n,m) \to {\L}^\prime(d+1,m_0,n,m)
\to \O_H(d+1) \to 0.
\]
Since $H^1(H,\O_H(d+1))=0$,
and by assumption
$H^1(({\bP}^\prime,{\L}^\prime(d,m_0,n,m))=0$,
we have the assertion.

Statement (c) is similarly proved, using the same argument
to reduce to $m_0'=m_0$, $n'=n$, and $m'=m$;
then the same exact sequence, at the $H^0$ level,
proves the result.

For (d) we note that $\v(d,0,2,d) = \v(2d,0,5,d) = d(1-d)/2 < 0$
but through $2$ general points there is always the line counted
with multiplicity $d$,
and through $5$ general points the conic counted with multiplicity $d$.
\end{pf}

Also, the case of large $m_0$ is easy to understand.
Denote by $T_i$ the line joining $p_0$ to $p_i$.

\begin{lemma}\label{basic2}\mbox{}
\begin{itemize}
\item[a.] $\L(d,m_0,n,m)$ is non-special if $m_0 > d$
($\l = \e = -1$ and the system is empty).
\item[b.] $\l(d,d,n,m) = \max\{-1,d-nm\}$ and
$\L(d,d,n,m)$ is special
if and only if $n\geq 1$, $m \geq 2$, and $d \geq nm$.
\item[c.] If $m+m_0 \geq d+1$ then
\[
\l(d,m_0,n,m) = \l(d-n(m+m_0-d),m_0-n(m+m_0-d),n,d-m_0).
\]
\item[d.] If $m \geq 2$ then $\l(d,d-1,n,m) = \max\{-1,2d-2nm+n\}$.
\end{itemize}
\end{lemma}

\begin{pf}
Statement (a) is obvious.
To see (c), note that $m+m_0 \geq d+1$ implies that
each line $T_i$ must be a component
of any divisor in the linear system;
therefore factoring out these $n$ lines implies that
$\l(d,m_0,n,m) = \l(d-n,m_0-n,n,m-1)$.
This is iterated $m+m_0-d$ times,
leading to
$\l(d,m_0,n,m) = \l(d-n(m+m_0-d),m_0-n(m+m_0-d),n,m-(m+m_0-d))$
as claimed.
Now (b) follows from (c), since if $m \geq 1$ then (c) gives
$\l(d,d,n,m) = \l(d-nm,d-nm,n,0)$.
The final statement is also proved by applying (c)
and using the Multiplicity One Lemma \ref{multone}.
\end{pf}

Later we will find it useful to apply Cremona transformations
to the plane to relate different linear systems
and to compute their dimensions.
Let $p_0,p_1,\ldots,p_n$ be general points in the plane,
and $C$ be a curve of degree $d$
having multiplicity $m_\ell$ at $p_\ell$ for each $\ell \geq 0$.
Fix three of the points $p_i$, $p_j$, and $p_k$;
note that $2d \geq m_i+m_j+m_k$ since otherwise the linear system
of such curves would have negative intersection number 
with every conic through these three points,
and therefore would be empty.
Then the effect of performing a quadratic Cremona transformation
based at these three points $p_i$, $p_j$, and $p_k$
is to transform $C$ to a curve of degree $2d-m_i-m_j-m_k$,
and having multiplicities at least
$d-m_j-m_k$ at $q_i$,
$d-m_i-m_k$ at $q_j$,
$d-m_i-m_j$ at $q_k$,
and $m_\ell$ at $p_\ell$ for $\ell \neq i,j,k$
(where $q_i$, $q_j$, and $q_k$
are the images of the three lines joining the
three points).

If we consider the entire linear system of such curves,
then it is clear that the dimension of the linear system
does not change upon performing the Cremona transformation.
Moreover, if the linear system contains irreducible curves
before applying the Cremona transformation,
it will contain irreducible curves after
applying the Cremona transformation,
except possibly the inclusion of one or more of the three
lines joining the three points as base locus.
In addition, the virtual dimension does not change
upon applying the Cremona transformation, if all of the numbers
involved are nonnegative.
Therefore Cremona transformation are a useful tool for analysing
the speciality of systems in certain situations.

\section{The Degeneration of the Plane}

In this section we describe the degeneration of the plane
which we use in the analysis.
It is related to that used by Ran \cite{ran}
in several enumerative applications.

Let $\Delta$ be a complex disc around the origin.
We consider the product $V = \bP^2 \times \Delta$,
with the two projections $p_1:V \to \Delta$ and $p_2:V \to \bP^2$.
We let $V_t = \bP^2\times \{t\}$.

Consider a line $L$ in the plane $V_0$ and blow it up
to obtain a new three-fold $X$ with maps
$f:X\to V$,
$\pi_1 = p_1 \circ f: X \to \Delta$,
and $\pi_2 = p_2 \circ f: X \to \bP^2$.
The map $\pi_1:X \to \Delta$
is a flat family of surfaces over $\Delta$.
We denote by $X_t$ the fibre of $\pi_1$ over $t \in \Delta$.
If $t \neq 0$, then $X_t = V_t$ is a plane $\bP^2$.
By contrast $X_0$ is the union
of the proper transform $\bP$ of $V_0$
and of the exceptional divisor $\bF$ of the blow-up.
It is clear that $\bP$ is a plane $\bP^2$
and $\bF$ is a Hirzebruch surface $\bF_1$,
abstractly isomorphic to the plane blown up at one point.
They are joined transversally along a curve $R$
which is a line $L$ in $\bP$
and is the exceptional divisor $E$ on $\bF$.

Notice that the Picard group of $X_0$
is the fibered product of $\Pic(\bP)$ and $\Pic(\bF)$
over $\Pic(R)$.
In other words, giving a line bundle $\mathcal{X}$ on $X_0$
is equivalent to giving a line bundle $\mathcal{X}_{\bP}$ on  ${\bP}$
and a line bundle $\mathcal{X}_{\bF}$ on  $\bF$
whose restrictions to $R$ agree.
The Picard group $\Pic(\bP)$ is generated by $\O(1)$,
while the Picard group of $\Pic(\bF)$ is generated by
the class $H$ of a line and the class $E$ of the exceptional divisor.
Since $H\cdot R = 0$ and $E\cdot R = -1$,
we have $\O_\bF(H)|_R \cong \O_R$ and $\O_\bF(E)|_R \cong \O_R(-1)$.
Hence in order that the restrictions to $R$ agree,
one must have
$\mathcal{X}_{\bP}\cong{\O}_{\bP}(d)$ and
$\mathcal{X}_{\bF}\cong{\O}_{\bF}(cH-dE)$ for some $c$ and $d$.
We will denote this line bundle on $X_0$ by $\mathcal{X}(c,c-d)$.

The normal bundle of $\bP$ in the $3$-fold $X$ is $-L$;
the normal bundle of $\bF$ in $X$ is $-E$.
Hence for example the bundle $\O_X(\bP)$
restricts to $\bP$ as ${\O}_{\bP}(-1)$
and restricts to $\bF$ as ${\O}_{\bF}(E)$.

Let ${\O}_X(d)$ be the line bundle $\pi_2^*({\O}_{\bP^2}(d))$.
If $t \not=0$, then the restriction of $\O_X(d)$ to $X_t\cong \bP^2$
is isomorphic to ${\O}_{\bP^2}(d)$
whereas the restriction of ${\O}_X(d)$ to $X_0$ is
the bundle $\mathcal{X}(d,0)$,
(whose restriction to $\bP$ is the bundle ${\O}_{\bP}(d)$
and whose restriction to $\bF$ is the bundle ${\O}_{\bF}(dH-dE)$).

Let us denote by ${\O}_X(d,k)$ the line bundle
${\O}_X(d)\otimes {\O}_X(k\bP)$.
The restriction of ${\O}_X(d,k)$ to $X_t$, $t\not=0$, is still the same,
i.e. it is isomorphic to ${\O}_{\bP^2}(d)$,
but the restriction to $X_0$ is now different:
it is isomorphic to  $\mathcal{X}(d,k)$
(whose restriction to $\bP$ is the bundle ${\O}_{\bP}(d-k)$
and whose restriction to $\bF$ is the bundle ${\O}_{\bF}(dH-(d-k)E)$).

We therefore see that all of the bundles $\mathcal{X}(d,k)$ on $X_0$
are flat limits of the bundles $\O_{\bP^2}(d)$
on the general fiber $X_t$ of this degeneration.

Fix a positive integer $n$
and another non-negative integer $b\leq n$.
Let us consider $n-b+1$ general points $p_0, p_1, \ldots, p_{n-b}$ in $\bP$
and $b$ general points $p_{n-b+1},...,p_n$ in $\bF$.
We can consider these points as limits
of $n$ general points $p_{0,t}, p_{1,t}, \ldots, p_{n,t}$ in $X_t$.
Consider then the linear system  ${\L}_t(d,m_0,n,m)$
which is the system ${\L}(d,m_0,n,m)$ in $X_t \cong {\bP}^2$
based at the points $p_{0,t}, p_{1,t},...,p_{n,t}$.

We now also consider the linear system
${\Lo}:={\Lo}(d,k,m_0,n,b,m)$ on $X_0$
which is formed by the divisors in $|\mathcal{X}(d,k)|$
having a point of multiplicity $m_0$ at $p_0$
and points of multiplicity $m$ at $p_1,...,p_n$.
According to the above considerations,
any one of the systems ${\Lo}(d,k,m_0,n,b,m)$
can be considered as a flat limit on $X_0$
of the system ${\L}_t(d,m_0,n,m)={\L}(d,m_0,n,m)$.
We will say that ${\Lo}$ is obtained from $\L$
by a \emph{$(k,b)$-degeneration}.

We note that the system $\Lo$ restricts to $\bP$
as a system $\LP$ of the form $\L(d-k,m_0,n-b,m)$
and $\Lo$ restricts to $\bF$
as a system $\LF$ of the form $\L(d,d-k,b,m)$.
Indeed, at the level of vector spaces,
the system $\Lo$ is the fibered product
of $\LP$ and $\LF$ over the restricted system on $R$,
which is $\O_R(d-k)$.
Specifically,
if $\LP$ is the projectivization of the vector space $W_\bP$,
and $\LF$ is the projectivization of the vector space $W_\bF$,
then by restriction to the double curve $R$ we have maps
$W_\bP \to H^0(R,\O(d-k))$ and $W_\bF \to H^0(R,\O(d-k))$,
and the fibered product $W = W_\bP \times_{H^0(R,\O(d-k))}W_\bF$
gives the linear system $\Lo = \bP(W)$ as its projectivization.

Since the linear system $\Lo$ is a linear system on a reducible scheme,
its elements come in three types.
The first type of element of $\Lo$
consists of a divisor $C_\bP$ on $\bP$ in the system $|(d-k)H|$
and a divisor $C_\bF$ on $\bF$ in the system $|dH-(d-k)E|$
(both of which satisfying the multiple point conditions)
which restrict to the same divisor on the double curve $R$.
We will then say that $C_\bP$ and $C_\bF$ \emph{match}
to give a divisor in $\Lo$.

The second type is a divisor corresponding to a section of the bundle
which is identically zero on $\bP$, and gives a general divisor in
the system $\L(d,d-k,b,m)$ on $\bF$
which contains the double curve $E$ as a component;
that is, an element of the system $E + \L(d,d-k+1,b,m)$.

The third type is the opposite,
corresponding to a section of the bundle
which is identically zero on $\bF$,
and gives a general divisor in
the system $\L(d-k,m_0,n-b,m)$ on $\bP$
which contains the double curve $L$ as a component;
that is, an element of the system $L + \L(d-k-1,m_0,n-b,m)$.

We denote by $\lo$ the dimension of the linear system $\Lo$ on $X_0$.
By semicontinuity, this dimension $\lo$
is at least that of the linear system on the general fiber,
i.e.,
\[
\lo = \dim(\Lo) \geq \l(d,m_0,n,m).
\]
Therefore we have the following:

\begin{lemma}
\label{lo=Ethenreg}
If $\lo = \e(d,m_0,n,m)$
then the system $\L(d,m_0,n,m)$ is non-special.
\end{lemma}

The basis of our method is to compute $\lo$ by a recursion.
The easy case is to compute this dimension
when all divisors in the linear system are of the second or third type,
that is, come from sections
which are identically zero on one of the components $\bP$ or $\bF$.
In this case one simply obtains the dimension of the linear system
on the other component, which gives us the following.

\begin{lemma}
\label{dimLo1}
Fix $d$, $k$, $m_0$, $n$, $b$, and $m$.
\begin{itemize}
\item[a.] If $\l(d-k,m_0,n-b,m) < 0$ then
$\lo = \l(d,d-k+1,b,m)$.
\item[b.] If $\l(d,d-k,b,m) < 0$ then
$\lo = \l(d-k-1,m_0,n-b,m)$.
\end{itemize}
\end{lemma}

We will define $\lohat$ to be the dimension
of the linear system $\Lohat$
of divisors in $\Lo$ which have the double curve $R$ as a component.

We need to extend Lemma \ref{dimLo1}
to handle the cases when there are divisors in $\Lo$
which are not identically zero on either component.
Fix $d$, $k$, $m_0$, $n$, $b$, and $m$.
We will refer to the system $\L = \L(d,m_0,n,m)$
as the \emph{general system}.
The system on $\bP$
restricts to a system $\RP$ on the double curve $R=L$,
and the kernel is, at the level of linear systems,
the system $\LPhat = \L(d-k-1,m_0,a,m)$.
Similarly, the system on $\bF$
restricts to a system $\RF$ on the double curve $R=E$,
and the kernel is, at the level of linear systems,
the system $\LFhat = \L(d,d-k+1,b,m)$.

We denote by
\[
\begin{array}{ll}
\v = \v(d,m_0,n,m) &
\text{ the virtual dimension of the general system,} \\
\vP = \v(d-k,m_0,n-b,m) &
\text{ the virtual dimension of the system on } \bP, \\
\vF = \v(d,d-k,b,m) &
\text{ the virtual dimension of the system on } \bF, \\
\vPhat = \v(d-k-1,m_0,n-b,m) &
\text{ the virtual dimension of the kernel system on } \bP, \\
\vFhat = \v(d,d-k+1,b,m) &
\text{ the virtual dimension of the kernel system on } \bF, \\
\l = \l(d,m_0,n,m) &
\text{ the dimension of the general system}, \\
\lP = \l(d-k,m_0,n-b,m) &
\text{ the dimension of the system on } \bP, \\
\lF = \l(d,d-k,b,m) &
\text{ the dimension of the system on } \bF, \\
\lPhat = \l(d-k-1,m_0,n-b,m) &
\text{ the dimension of the kernel system on } \bP, \\
\lFhat = \l(d,d-k+1,b,m) &
\text{ the dimension of the kernel system on } \bF, \\
\rP = \lP-\lPhat-1 &
\text{ the dimension of the restricted system $\RP$}\\
&\text{ on $R=L$ from $\bP$, and} \\
\rF = \lF-\lFhat-1 &
\text{ the dimension of the restricted system $\RF$}\\
&\text{ on $R=E$ from $\bF$}. \\
\end{array}
\]

One has the following lemma,
whose immediate proof can be left to the reader:

\begin{lemma}
\label{identities}
The following identities hold:
\begin{itemize}
\item[a.] ${\vP}+{\vF}=\v + d-k$.
\item[b.] ${\vPhat}+{\vF}={\v}-1$.
\item[c.] ${\vP}+{\vFhat}={\v}-1$.
\item[d.] ${\lohat}={\lPhat}+{\lFhat}+1$.
\end{itemize}
\end{lemma}

The restricted systems $\RP$ and $\RF$ on the double curve
may intersect in various dimensions a priori.
The dimension $\lo$ of the linear system $\Lo$ on $X_0$
depends on the dimension of this intersection,
since $\Lo$ is obtained as a fibered product.
Returning to the notation above,
we have that $\Lo = \bP(W)$ where $W$ is the fibered product
$W_\bP \times_{H^0(R,\O(d-k))} W_\bF$.
Hence at the level of vector spaces
\[
W = \{(\alpha,\beta) \in W_\bP \times W_\bF \;|\;
\alpha|_R = \beta|_R\}.
\]
Let $W_R$ be the vector space corresponding to the intersection
of the restricted systems $\RP\cap\RF$, so that $\RP\cap\RF = \bP(W_R)$.
Then an element of $W$ is determined by first choosing an element
$\gamma \in W_R$,
then choosing pre-images $\alpha \in W_\bP$ and $\beta \in W_\bF$
of $\gamma$.
Using vector space dimensions, the choice of $\gamma$
depends on $1+\dim(\RP\cap\RF)$ parameters,
and then once $\gamma$ is chosen the choice of $\alpha$
depends on the vector space dimension of the kernel system,
which is $1+\lPhat$, and similarly the choice of $\beta$
depends on $1+\lFhat$ parameters.
Therefore $\dim(W) = 1+\dim(\RP\cap\RF) + 1+\lPhat + 1+\lFhat$.
Projectivizing gives the dimension of the linear system $\Lo$,
and we have proved the following:

\begin{lemma}
\label{dimLo2}
With the above notation,
\[
\lo = \dim(\RP\cap\RF) + \lPhat + \lFhat + 2.
\]
\end{lemma}

This is the extension of Lemma \ref{dimLo1} which we were seeking.


\section{The Transversality of the Restricted Systems}

It is clear from the previous Lemma
that the computation of $\lo$
depends on the knowledge of the dimension
of the intersection $\RP\cap\RF$
of the restricted linear systems.
The easiest case to handle would be if these two systems
were transverse (as linear subspaces of the projective space
of divisors of degree $d-k$ on the double curve $R$);
then a formula for the dimension of the intersection is immediate.
It turns out that this is always the case,
which is a consequence of the following Proposition,
first proved to our knowledge by Hirschowitz in
\cite{hirschowitz2}, using the Borel fixed point theorem.
Our proof is a variation on the theme using
the finiteness of inflection points of linear systems.

\begin{proposition}
\label{SL2prop}
Let $G = PGL(2,\mathbb{C})$ be the automorphism group of $\mathbb{P}^1$.
Let $X$ be the linear system of divisors of degree $d$ on $\mathbb{P}^1$.
Note that $G$ acts naturally on $X$,
and on linear subspaces of $X$ of any dimension.
Then for any two nontrivial linear subspaces $V$ and $W$ of $X$,
there is an element $g \in G$ such that $V$ meets $gW$ properly.
\end{proposition}

\begin{pf}
It suffices to prove the assertion when $V$ and $W$ have complementary
dimensions $k$ and $d-k-1$ respectively.
We argue in this case by contradiction:
suppose that for every $g \in G$, the intersection $V\cap gW$ is nonempty.
Fix a general coordinate system $[x:y]$ on $\mathbb{P}^1$,
and consider the element $g_t \in G$ given by
$g_t[x:y] = [tx:t^{-1}y]$.
Suppose that in this coordinate system
we have a basis $\{v_0,\dots,v_k\}$ for $V$
and a basis $\{w_{k+1},\ldots,w_d\}$ for $W$.
Write each $v_i$ as a polynomial in $x$ and $y$ as
$v_i=\sum_j a_{ij}x^jy^{d-j}$
and similarly write each $w_i$ as $w_i=\sum_j b_{ij}x^jy^{d-j}$.
Note that with this notation we have that the subspace $g_tW$
has as its basis the polynomials
$g_tw_i = \sum_j t^{2j-d}b_{ij}x^jy^{d-j}$.

Let $A$ be the $(k+1) \times (d+1)$ matrix of the $a_{ij}$ coefficients,
and let $B_t$ be the $(d-k)\times (d+1)$ matrix
of the $t^{2j-d}b_{ij}$ coefficients.
Notice that $B:=B_1$ is the matrix of coefficients
for the original subspace $W$.
Let $C_t$ be the square matrix whose with $A$ as its first $k+1$ rows
and $B_t$ as its last $d-k$ rows.

Since the subspaces $V$ and $g_tW$ intersect nontrivially,
they cannot span the whole space $X$;
hence the matrix of coefficients $C_t$
must have trivial determinant.
Note that this determinant is a Laurent polynomial in $t$,
and hence each coefficient of $t$ in this polynomial must vanish.

By expressing this determinant in the Laplace expansion
using the minors of the first $k+1$ rows
against the minors of the last $d-k$ rows,
we see that the top coefficient of the determinant
is the product of the minor of $A$ using the first $k+1$ columns
with the minor of $B$ using the last $d-k$ columns.
Since this coefficient is zero, we have that
either the first minor of $A$ is zero
or the last minor of $B$ is zero.

If the first minor of $A$ is zero,
then there exists in $V$ a polynomial
whose first $k+1$ coefficients are zero,
and hence vanishes at $[0:1]$ to order at least $k+1$.
Hence $[0:1]$ would be an inflection point for the system $V$.

Similarly, if the last minor of $B$ is zero,
we conclude in the same way
that the point $[1:0]$ is an inflection point
for the system $W$.

Since the coordinate system was chosen to be general,
we see that there are infinitely many inflection points
for at least one of the two systems.
This is a contradiction, finishing the proof.
\end{pf}

Note that the given any automorphism $g$ of a line in the plane,
there is a lift of $g$ to an automorphism of the plane fixing the line.
By using this and the previous Proposition,
one immediately deduces the following.

\begin{corollary}
The restricted systems $\RP$ and $\RF$ on the double line
intersect properly.
\end{corollary}

The previous Corollary, combined with Lemma \ref{dimLo2},
gives the formula for $\lo$:

\begin{proposition}
\label{dimLo3}
\mbox{}
\begin{itemize}
\item[(a)] If $\rP+\rF \leq d-k-1$, then
\[
\lo = \lPhat + \lFhat + 1.
\]
\item[(b)] If $\rP+\rF \geq d-k-1$, then
\[
\lo = \lP+ \lF - d + k.
\]
\end{itemize}
\end{proposition}

\begin{pf}
If $\rP+\rF \leq d-k-1$, then the transversality
of $\RP$ and $\RF$ implies that $\RP\cap\RF$ is empty, of dimension $-1$.
This gives (a), using Lemma \ref{dimLo2}.
If $\rP+\rF \geq d-k-1$, then $\dim(\RP\cap\RF) = \rP+\rF-d+k$
(again using the transversality)
and by Lemma \ref{dimLo2}, we have

\begin{align*}
\lo &= \rP+\rF-d+k  +  \lPhat  + \lFhat + 2\\
&= \lP +\lF - d + k
\end{align*}
using the definition of $\rP=\lP-\lPhat-1$ and $\rF=\lF-\lFhat-1$.
\end{pf}

Our method for proving that the system $\L$ is non-special
is to find appropriate integers $k$, $a$, and $b$ with $n = a+b$
such that $\lo = \e(d,m_0,a+b,m)$,
and to invoke Lemma \ref{lo=Ethenreg}.
The computation of $\lo$ is done
by recursively using Proposition \ref{dimLo3}.

The dimension computed in part (b) of the Proposition
is, miraculously,
the virtual dimension of the system on the general fiber,
if each of the systems involved in (b) has the virtual dimension,
using Lemma \ref{identities}(a).
Therefore (b) is useful for proving that $\L$
has the correct minimal dimension and is therefore non-special.
Statement (a) is more useful for proving that $\L$ is empty.

The following makes these two strategies explicit.

\begin{corollary}
\label{cordimLk}
Fix $d$, $m_0$, $n$, and $m$.
\begin{itemize}
\item[a.] Suppose that positive integers $k$ and $b$ exist,
with $0 < k < d$ and $0 < b < n$, such that
\begin{enumerate}
\item $\L(d-k-1,m_0,n-b,m)$ is empty,
\item $\L(d,d-k+1,b,m)$ is empty, and
\item $\dim \L(d-k,m_0,n-b,m) + \dim \L(d,d-k,b,m) \leq d-k-1$.
(This is automatic if $\v(d,m_0,n,m) \leq -1$
and both these systems are non-special
with virtual dimension at least $-1$.)
\end{enumerate}
Then $\L(d,m_0,n,m)$ is empty (and therefore non-special).
\item[b.] Suppose that positive integers $k$ and $b$ exist,
with $0 < k < d$ and $0 < b < n$, such that
\begin{enumerate}
\item $\v(d-k,m_0,n-b,m) \geq -1$ and $\L(d-k,m_0,n-b,m)$ is non-special,
\item $\v(d,d-k,b,m) \geq -1$ and $\L(d,d-k,b,m)$ is non-special, and
\item $\dim \L(d-k-1,m_0,n-b,m) + \dim \L(d,d-k+1,b,m)
\leq \v(d,m_0,n,m)-1$.
(This is automatic if both these systems are non-special
with virtual dimension at least $-1$.)
\end{enumerate}
Then $\L(d,m_0,n,m)$ is non-special,
with virtual dimension $\v(d,m_0,n,m)$ at least $-1$.
\end{itemize}
\end{corollary}

\begin{pf}
In case (a), the first two hypotheses say that
$\lPhat = \lFhat = -1$,
so that $\rP = \lP$  and $\rF = \lF$.
Hence by the third hypothesis
\[
\rP+\rF = \lP + \lF \leq d-k-1
\]
so that using Proposition \ref{dimLo3}(a)
we have $\lo = \lPhat + \lFhat + 1 = -1$,
proving that $\L$ is empty by semicontinuity.
The parenthetical statement in the third hypothesis
follows from Lemma \ref{identities}(a).

For (b), the first two hypotheses say that $\lP=\vP$ and $\lF=\vF$;
then 
\begin{align*}
\rP+\rF &= \lP - \lPhat + \lF - \lFhat - 2 \\
&= \vP - \lPhat + \vF - \lFhat - 2 \\
&\geq \vP+\vF - 1 - \v \;\;\;\text{ using the third hypothesis}\\
&= d-k-1
\end{align*}
using Lemma \ref{identities}(a).
Hence Proposition \ref{dimLo3}(b) implies that
$\lo = \lP+\lF-d+k = \vP + \vF - d + k = \v(d,m_0,n,m)$
again using Lemma \ref{identities}(a).
The parenthetical statement in the third hypothesis
follows from Lemma \ref{identities}(a), (b), and (c).
\end{pf}

\section{$(-1)$-Special Systems and the Main Conjecture}

A linear system $\L(d,m_0,n,m)$
with $\L^2 = -1$ and $g_\L = 0$
will be called a \emph{quasi-homogeneous $(-1)$-class}.
By (\ref{RR}), we see that $\v = 0$,
so that every quasi-homogeneous $(-1)$-class is effective.

Suppose that $A$ is an irreducible rational curve
and is a member of a linear system $\L=\L(d,m_0,n,m)$,
and suppose that on the blowup $\bP^\prime$ of the plane
the proper transform of $A$ is smooth, of self-intersection $-1$.
We say then that $A$ is a \emph{$(-1)$-curve}.
In this case $\L$ is a quasi-homogeneous $(-1)$-class.
Indeed, if this happens, then $\L = \{A\}$:
if $D \in \L$, then $D\cdot A < 0$ on $\bP^\prime$,
so that $D$ must contain $A$ as a component,
and then be equal to $A$ since they have the same divisor class.
Therefore such a linear system $\L$ is non-special, of dimension $0$.
A quasi-homogeneous $(-1)$-class
containing a $(-1)$-curve
will be called an \emph{irreducible $(-1)$-class}.

Let $A$ be a $(-1)$-curve and suppose that
$2A$ is a member of a linear system $\L$.
Then $\L^2 = -4$ and $g_\L = -2$
so by (\ref{RR}) $\v = -1$ and the system is expected to be empty;
however it clearly contains the divisor $2A$
(and is equal in fact to $\{2A\}$).
Therefore such a linear system is special.

More generally we have the following.

\begin{lemma}
\label{minusonespecial}
Let $\L$ be a nonempty linear system,
Suppose that $A_1,\ldots,A_r$ are $(-1)$-curves which meet $\L$ negatively;
write $\L \cdot A_j = - N_j$ with $N_j \geq 1$ for each $j$.
Then:
\begin{itemize}
\item[(a)] $\L$ contains $\sum_j N_jA_j$ as a fixed divisor.
\item[(b)] If $i \neq j$ then $A_i\cdot A_j = 0$.
\item[(c)] If $\M = \L - \sum_j N_j A_j$ is the residual system, then
\[
\v(\M) - \v(\L) = \sum_j N_j(N_j-1)/2.
\]
\item[(d)] If $N_j \geq 2$ for some $j$ then $\L$ is special.
\end{itemize}
\end{lemma}

\begin{pf}
We work on the blowup $\bP^\prime$,
and we note that since $\L\cdot A_j = -N_j$,
$N_jA_j$ is certainly a fixed divisor;
hence we have $\L = \sum_j N_jA_j + \M$
for some linear system $\M$, with $\dim(\L) = \dim(\M)$,
which proves (a).

If a linear system $\L$ meets two $(-1)$-curves $A$ and $B$
negatively,
then both $A$ and $B$ must be part of the fixed divisor
of $\L$.
If $A$ and $B$ meet then since
$\v(A+B)=A\cdot B$, we would have that $A+B$ moves,
and could not be part of the fixed divisor of $\L$.
This proves (b).

Note also that $\M^2 = \L^2 + \sum_j N_j^2$,
and $\M\cdot K = \L \cdot K + \sum_j N_j$.
Then
\begin{align*}
\v(\M) - v(\L) &= (\M^2-\M\cdot K)/2 - (\L^2-\L\cdot K)/2 \\
&= (\L^2+\sum_j N_j^2-\L\cdot K-\sum_j N_j)/2 - (\L^2- \L \cdot K)/2 \\
&= \sum_j N_j(N_j-1)/2 > 0
\end{align*}
which proves (c).
If any $N_j \geq 2$ then we see that $v(\M) > \v(L)$ and hence
\[
\dim(\L) = \dim(\M) \geq \v(\M) > \v(\L)
\]
proving the speciality of $\L$.
\end{pf}

The above Lemma suggests the following.

\begin{definition}
A linear system $\L$ is \emph{$(-1)$-special}
if there are $(-1)$-curves $A_1,\ldots,A_r$
such that $\L\cdot A_j = -N_j$ with $N_j \geq 1$ for every $j$
and $N_j \geq 2$ for some $j$,
with the residual system $\M = \L - \sum_j N_j A_j$
having non-negative virtual dimension $\v(\M) \geq 0$,
and having non-negative intersection with every $(-1)$-curve.
\end{definition}

We remark that if $\L$ is a linear system
satisfying all of the hypotheses of the above definition
except the last one,
then in fact it is $(-1)$-special;
the residual system can meet only finitely many
additional $(-1)$-curves negatively,
and after adding these to the fixed part
we get a residual system satisfying the final condition.

By the Lemma, every $(-1)$-special system is special;
the condition that $\v(\M) \geq 0$ implies that $\M$
is non-empty, and hence that $\L$ is non-empty.

The following conjecture is a restatement of a conjecture
of Hirschowitz (see \cite{hirschowitz3}),
also made by Harbourne
(see \cite{harbourne1} and \cite{harbourne2}).

\begin{mainconj}
Every special system is $(-1)$-special.
\end{mainconj}

\section{Quasi-homogeneous $(-1)$-classes}

It is clear from the previous section
that a classification of $(-1)$-classes is important
in understanding speciality of linear systems.
Fortunately it is not hard
to classify all quasi-homogeneous $(-1)$-classes,
as we now do.

Suppose that $\L(d,m_0,n,m)$ is a quasi-homogeneous $(-1)$-class.
Then
\begin{equation}
\label{QH1}
d^2-m_0^2-nm^2 = -1
\end{equation}
(since $\L^2 = -1$)
and
\begin{equation}
\label{QH2}
3d-m_0-nm=1
\end{equation}
(which is equivalent to the genus condition,
and is specifically the condition that $\L\cdot K = -1$
on the blowup surface).
Solving (\ref{QH2}) for $m_0$ gives $m_0=3d-nm-1$
and plugging this into (\ref{QH1}) yields
\begin{equation*}
8d^2 - 6dnm + n^2m^2 - 6d + 2nm + nm^2 = 0
\end{equation*}
which can be rewritten as
\begin{equation}
\label{QH4}
(4d-nm)(2d-nm) + (m-1)(4d-nm) - (2m+1)(2d-nm) = 0.
\end{equation}
This suggests the change of variables
\[
u = 4d - nm, \;\;\; v = 2d - nm
\]
so that (\ref{QH4}) now becomes
\begin{equation}
\label{QH5}
uv + (m-1)u - (2m+1)v = 0.
\end{equation}
Reversing this change of coordinates gives
\begin{equation}
\label{QH7}
d = \frac{u-v}{2},\;\;\;
nm = u-2v,\;\;\;
m_0 = \frac{u+v}{2} - 1.
\end{equation}
Hence we seek integral solutions $(u,v)$ to the equation (\ref{QH5})
with $u\equiv v \mod{2}$ (so that $d$ and $m_0$ are integers)
and $m|(u-2v)$ (so that $n$ is an integer),
and all quantities $u-v$, $u-2v$, and $u+v$ positive.

If $m=1$, the curve (\ref{QH5}) is $uv-3v=0$,
so either $u=3$ or $v=0$.
If $u=3$ then the positivity conditions are that
$3-v>0$, $3-2v>0$, and $3+v>0$, so $-3 < v < 3/2$ and must be odd;
only $v = \pm 1$ are possibilities.
The solution $(3,-1)$ gives $(d,m_0,n,m)=(2,0,5,1)$,
and the solution $(3,1)$ gives $(d,m_0,n,m)=(1,1,1,1)$.
If $v=0$ we only must have $u>0$ and even, say $u=2e$.
This gives $(d,m_0,n,m)=(e,e-1,2e,1)$, for any $e \geq 1$.
These are all the solutions with $m=1$.

{}From now on we assume that $m \geq 2$.
In this case the hyperbola (\ref{QH5}) in the $(u,v)$ plane
has the horizontal asymptote $v = 1-m$,
the vertical asymptote $u = 2m+1$,
and passes through the origin with slope $(m-1)/(2m+1) < 1/2$.
This slope condition implies that in the third quadrant
the hyperbola lies entirely above the line $v = u/2$.

Now $nm = u-2v > 0$ so $v < u/2$.
Moreover $m_0 \leq d-1$ hence $(u+v)/2 -1 \leq (u-v)/2 -1$,
implying that $v \leq 0$.
These two inequalities imply that the only integral points of interest
lie on the branch of the hyperbola in the fourth quadrant,
with $v < 1-m$ and $u > 2m+1$.  Hence we may assume that
\begin{equation}
\label{QH8}
u \geq 2m+2
\;\;\;\mathrm{and}\;\;\;
v \leq -m.
\end{equation}

Finally make the change of coordinates
\begin{equation}
x = u-2m-1, \;\;\; y = 1-m-v, \;\;\; u=x+2m+1,\;\;\; v = 1-m-y
\end{equation}
which transforms the hyperbola (\ref{QH5}) into
\begin{equation}
\label{QH9}
xy = (m-1)(2m+1).
\end{equation}
The branch of (\ref{QH9}) corresponding to the branch of (\ref{QH5})
in the fourth quadrant is the one in the first quadrant,
with $x\geq 1$ and $y\geq 1$.
Clearly the integral points on (\ref{QH9})
come simply from the possible factorizations of $(m-1)(2m+1)$.

This gives the following classification.

\begin{proposition}
\label{QH-1general}
The quasi-homogeneous $(-1)$-classes are the classes $\L(d,m_0,n,m)$
with $(d,m_0,n,m)$ on the following list:
\begin{itemize}
\item[(a)] $(2,0,5,1)$ and $(1,1,1,1)$.
\item[(b)] $(e,e-1,2e,1)$ with $e \geq 1$.
\item[(c)] For any $m \geq 2$, and any $x\geq 1$, $y\geq 1$ with
\begin{itemize}
\item[(i)]  $xy = (m-1)(2m+1)$,
\item[(ii)] $x+m \geq y$,
\item[(iii)] $x-y \equiv m \mod{2}$, and
\item[(iv)] $m | x+2y-1$,
\end{itemize}
the four-tuple
\[
(\frac{x+y+3m}{2}, \frac{x-y+m}{2}, \frac{x+2y-1}{m}+4, m).
\]
\end{itemize}
\end{proposition}

In part (c), condition (ii) is the non-negativity of $m_0$,
while conditions (iii) and (iv) are needed to insure
that $d$, $m_0$, and $n$ are integral.

It is easy to classify all homogeneous $(-1)$-classes
from the above Proposition.

\begin{corollary}
The classes
$\L(2,0,5,1)=\L(2,1,4,1)$ and $\L(1,0,2,1)=\L(1,1,1,1)$
are the only homogeneous $(-1)$-classes.
\end{corollary}

\begin{pf}
Clearly these are the only ones with $m=1$.
For $m \geq 2$, we must have the factors $x$ and $y$
satisfying $y=x+m$, and so $x(x+m)=(m-1)(2m+1)$.
If $x \leq m-1$ then the other factor $x+m$ must be at least $2m+1$,
a contradiction.
If $x \geq m+1$ then $x(x+m) \geq (m+1)(2m+1)$, which is too big.
Hence only $x=m$ is a possibility, which in fact does not work.
\end{pf}

\begin{example}
\label{QH-1extremal}
For any $m\geq 2$, set $x = (m-1)(2m+1)$ and $y = 1$.
Conditions (i), (ii), and (iii) clearly hold,
and $x+2y-1 = (2m^2 -m - 1) +2-1 = m(2m-1)$
so that also (iv) holds.
This gives
\[
d=m^2+m, \;\;\; m_0 = m^2-1,\;\;\; n=2m+3.
\]
\end{example}

\begin{example}
Fix an integer $z\neq 0$, let $m=4z^2+3z$,
$x = 8z^2+2z-1$ and $y = 4z^2+5z+1$.
This gives the solution
\[
d=12z^2+8z, \;\;\; m_0 = 4z^2-1,\;\;\; n=8,\;\;\; m=4z^2+3z.
\]
\end{example}


It is an exercise to check that the previous two examples
produce all of the quasi-homogeneous $(-1)$-classes with $n=8$.

\begin{example}
\label{QH-1list}
The following is a complete list
of all quasi-homogeneous $(-1)$-classes
with $m \leq 7$:

\begin{center}
\begin{tabular}{cccccc}
$d$ & $m_0$ & $n$ & $m$ & ($x$ & $y$) \\ \hline
1 & 1 & 1 & 1 & - & - \\
2 & 0 & 5 & 1 & - & - \\
$e\geq 1$ & $e-1$ & $2e$ & $1$ & - & - \\
6 & 3 & 7 & 2 & (5 & 1)\\
12 & 8 & 9 & 3 & (14 & 1)\\
20 & 15 & 11 & 4 & (27 & 1) \\
30 & 24 & 13 & 5 & (44 & 1) \\
42 & 35 & 15 & 6 & (65 & 1) \\
20 & 3 & 8 & 7 & (9 & 10) \\
27 & 17 & 9 & 7 & (30 & 3) 
\end{tabular}
\end{center}
\end{example}

We note that Cremona transformations may be used
to give a numerical criterion for deciding
when a quasi-homogeneous $(-1)$-class is an irreducible class:
if it can be transformed,
by a series of quadratic Cremona transformations,
to the class of a line through two points $\L(1,0,2,1)$.

\begin{proposition}\mbox{}
\begin{itemize}
\item[(a)] All quasi-homogeneous $(-1)$-classes
having $m =1$ are irreducible.
\item[(b)] All quasi-homogeneous $(-1)$-classes with $m \geq 2$
of the form $\L(d=m^2+m, m_0=m^2-1, n=2m+3, m)$
(obtained by the factorization $x = (m-1)(2m+1)$, $y=1$,
as in Example \ref{QH-1extremal}) are irreducible.
\item[(c)] All quasi-homogeneous $(-1)$-classes
having $m \leq 6$ are irreducible.
\item[(d)] The $(-1)$-class $\L(27,17,9,7)$ is not irreducible.
\end{itemize}
\end{proposition}

\begin{pf}
It is obvious that $\L(1,1,1,1)$ and $\L(2,0,5,1)$,
corresponding to a line through two points and a conic through $5$,
are irreducible.
To see that $\L(e,e-1,2e,1)$ is irreducible,
note that it is irreducible for $e=1$:
again this is a line through $2$ points.
For $e \geq 2$,
applying a Cremona transformation to $\L(e,e-1,2e,1)$
at the points $p_0$, $p_1$, $p_2$
transforms the system to $\L(e-1,e-2,2e-2,1)$,
and so by induction all these systems are irreducible.
This proves (a).

To prove (b), apply the quadratic Cremona transformation
to $\L(m^2+m, m^2-1, 2m+3, m)$ exactly $m+1$ times,
at $p_0, p_{2i-1}, p_{2i}$, for $i=1,\ldots m+1$.
It is easy to see that this transforms the system to
$\L(m+1,m,2m+2,1)$, which is irreducible by (a).

Part (c) now is a consequence of (a) and (b),
given the list of Example \ref{QH-1list}.

To prove (d), note that $\L(12,8,9,3)$ is irreducible,
and $\L(27,17,9,7)\cdot \L(12,8,9,3)
= 12\cdot 27 - 8\cdot 17 - 9\cdot 3 \cdot 7 = -1$,
so that if $A$ is the $(-1)$-curve in $\L(12,8,9,3)$,
then $A$ is a fixed curve of $\L(27,17,9,7)$.
The residual system is $\L(15,9,9,4)$,
which has virtual dimension $0$ and is therefore non-empty.
\end{pf}

Recall that we are interested in $(-1)$-curves
because they are useful in constructing special linear systems.
Suppose a quasi-homogeneous system $\L(d,m_0,n,m)$
meets negatively a $(-1)$-curve $A$
of degree $\delta$, having multiplicities
$\mu_0,\mu_1,\ldots,\mu_n$ at the points $p_0,\ldots,p_n$.
Since the points are general, a monodromy argument
implies that for any permutation $\sigma \in \Sigma_n$,
$\L$ also meets negatively
the $(-1)$-curve $A_\sigma$ of degree $\delta$,
having multiplicity $\mu_0$ at $p_0$,
and having multiplicities $\mu_{\sigma(i)}$ at $p_i$
for each $i \geq 1$.
There may of course be repetitions among the $A_\sigma$'s.

If this happens, no two of the $(-1)$-curves $A_\sigma$ can meet,
by Lemma \ref{minusonespecial}.
Hence, if $m+0\neq0$,
the Picard group of the blowup surface $\bP^\prime$
has rank $n+2$,
and there can be at most $n+1$ of these disjoint $(-1)$-curves.
In the homogeneous case where $m_0=0$,
we do not blow up $p_0$,
and the rank of the Picard group is only $n+1$;
therefore in this case there can be at most $n$ of these disjoint
$(-1)$-curves.

The sum of all of these $A_\sigma$'s must also be quasi-homogeneous,
and if there are $k$ of them,
is therefore of the form $\L(k\delta,k\mu_0,n,\mu^\prime)$
for some $\mu^\prime$.
An elementary counting argument shows that if the $\mu_i$'s
(for $i\geq 1$) occur in subsets of size $k_1\leq k_2\leq\dots\leq k_s$,
constant in each subset,
then the number of distinct $A_\sigma$'s is
\[
\frac{n!}{k_1!k_2!\cdots k_s!}.
\]
The only way this can be less than or equal to $n+1$
is if $s=1$ (and $A$ is then quasi-homogeneous)
or if $s=2$ and $k_1=1$, $k_2=n-1$.

The classification in case $s=1$ we have discussed above.
In the case $s=2$, there are exactly $n$ $A_\sigma$'s,
and the sum of the $A_\sigma$'s
is quasi-homogeneous, and is of the form
$\L(n\delta,n\mu_0,n,\mu_1+(n-1)\mu_2)$
if $\mu_i = \mu_2$ for $i \geq 2$.
The condition that $A$ and $A_\sigma$ do not meet
is that
\[
\delta^2 - \mu_0^2 - 2\mu_1\mu_2 - (n-2)\mu_2^2 = 0
\]
while that fact that $A$ is a $(-1)$-curve implies that
\[
\delta^2 - \mu_0^2 - \mu_1^2 - (n-1)\mu_2^2 = -1.
\]
Subtracting these two equations gives
\[
{(\mu_1-\mu_2)}^2 = 1
\]
so that $\mu_1 = \mu_2 \pm 1$.

We call such a system a \emph{quasi-homogeneous $(-1)$-configuration}.
We say that the configuration is \emph{compound}
if it consists of more than one $(-1)$-curve.

It is not hard to completely classify these classes for low $m$:

\begin{example}
\label{compoundQH-1list}
The following is a complete list
of all of the quasi-homogeneous $(-1)$-configurations $\L(d,m_0,n,m)$
with $m \leq 10$:

\begin{center}
\begin{tabular}{cccccccl}
$d$ & $m_0$ & $n$ & $m$ &
($\delta$ & $\mu_0$ & $\mu_1$ & $\mu_2=\cdots = \mu_n$) \\ \hline
$e \geq 2$ & $e$ & $e$ & 1 & (1 & 1 & 1 & 0) \\
3        & 0 & 3 & 2 & (1 & 0 & 0 & 1) \\
10       & 5 & 5 & 4 & (2 & 1 & 0 & 1) \\
12       & 0 & 6 & 5 & (2 & 0 & 0 & 1) \\
21      & 14 & 7 & 6 & (3 & 2 & 0 & 1) \\
18       & 6 & 6 & 7 & (3 & 1 & 2 & 1) \\
21       & 0 & 7 & 8 & (3 & 0 & 2 & 1) \\
36      & 27 & 9 & 8 & (4 & 3 & 0 & 1) \\
55      & 44 & 11 & 10 & (5 & 4 & 0 & 1) \\
\end{tabular}
\end{center}
\end{example}

Using the classification of quasi-homogeneous $(-1)$-classes,
one can easily give a complete classification of
homogeneous $(-1)$-configurations, i.e., those with $\mu_0=0$.

\begin{proposition}
\label{homog-1}
The following is a complete list of homogeneous $(-1)$-configurations:
\[
\begin{array}{cl}
\L(1,0,2,1) & \mathrm{\;\;\;not\;\;compound} \\
\L(2,0,5,1) & \mathrm{\;\;\;not\;\;compound} \\
\L(3,0,3,2) & \mathrm{\;\;\;compound\;\;with\;\;\;} \delta=1, n=3, \mu_1=0,
\mu_2=1\\
\L(12,0,6,5) & \mathrm{\;\;\;compound\;\;with\;\;\;} \delta=2, n=6, \mu_1=0,
\mu_2=1\\
\L(21,0,7,8) & \mathrm{\;\;\;compound\;\;with\;\;\;} \delta=3, n=7, \mu_1=2,
\mu_2=1\\
\L(48,0,8,17) & \mathrm{\;\;\;compound\;\;with\;\;\;} \delta=6, n=8, \mu_1=3,
\mu_2=2.
\end{array}
\]
\end{proposition}

\begin{pf}
The non-compound $(-1)$-curves $\L(1,0,2,1)$ and $\L(2,0,5,1)$
we have seen before.
If $A$ is a $(-1)$-curve producing a homogeneous $(-1)$-configuration,
then, with the notation above, $A$ has degree $\delta$,
one point of multiplicity $\mu_1$, and $n-1$ points of multiplicity
$\mu_2 = \mu_1 \pm 1$.
Therefore, shifting the indices of the points,
we see that $A$ is quasi-homogeneous, in the class
$\L(\delta,\mu_1,n-1,\mu_2)$.
Hence we may appeal to Proposition \ref{QH-1general}.

If $\mu_2 = 1$, then either $\mu_1=0$
(giving the two possibilities
$\delta = 1, n-1=2$, and the compound configuration $\L(3,0,3,2)$, or
$\delta=2$, $n-1=5$, and the compound configuration $\L(12,0,6,5)$)
or $\mu_1=2$
(giving $\delta=3$, $n-1=6$, and the compound configuration $\L(21,0,7,8)$).

If $\mu_2 \geq 2$, then the class $A$ comes from a factorization
$xy = (\mu_2-1)(2\mu_2+1)$,
and then $\mu_1 = (x-y+\mu_2)/2$.
This is $\mu_2\pm 1$ if and only if $x-y = \mu_2 \pm 2$.

Now $x = 2\mu_2+1$, $y = \mu_2-1$ is a factorization with $x-y = \mu_2+2$,
but now we look at the requirement
that $\mu_2$ divides $x+2y-1 = 4\mu_2-2$.
This forces $\mu_2 = 2$, giving $x=5$, $y=1$, and $\delta = 6, n-1=7$,
leading to the compound configuration $\L(48,0,8,17)$.

This is the only solution with $x-y=\mu_2+2$.
Hence what is left is to discuss the possible cases
with $x-y = \mu_2-2$.
To obtain such a factorization $x$, $y$,
we must have $x < 2\mu_2+1$ and $y > \mu_2-1$.
Neither $x = 2\mu_2$ nor $y = \mu_2$ are possible factors,
so in fact we must have $x \leq 2\mu_2-1$ and $y \geq \mu_2+1$.
However the difference $x-y$ being $\mu_2-2$
then forces $x = 2\mu_2-1$ and $y = \mu_2+1$,
whose product never equals $(\mu_2-1)(2\mu_2+1)$.
Thus there are no more homogeneous $(-1)$-configurations.
\end{pf}

\section{The Dimension for Large $m_0$}
\label{sec10}

In our application we will usually
make a $(k,b)$-degeneration with $k$ near $m$.
This leads to linear systems on $\bF$
which have the form $\L(d,m_0,b,m)$
with $m_0$ near $d-m$.
These linear systems may usually be effectively analyzed
by applying Cremona transformations,
since one of the multiplicities is so large with respect to the degree.

Let us first consider quasi-homogeneous systems of the form $\L(d,d-m,n,m)$.

\begin{lemma}
\label{m0=d-m_algorithm}
Fix $d \geq m\geq 0$.
Consider the general linear system $\L=\L(d,d-m,n,m)$,
whose dimension is $\l$.
\begin{itemize}
\item[(a)] If $m=0$ then $\L$ is non-special
and $\l=\v(d,d)=d$.
\item[(b)] If $m=1$ then $\L$ is non-special
and $\l = \e(d,d-1,n,1)=\max\{-1,2d-n\}$.
\item[(c)] If $n=0$ then $\L$ is non-special
and $\l = \v(d,d-m) = d+dm-m^2/2+m/2$.
\item[(d)] If $n=1$ then $\L$ is non-special
and $\l = \v(d,d-m,1,m) = d+m(d-m)$.
\item[(e)] If $n = 2$ and $m \leq d \leq 2m$ then
$\l = (d-m)(d-m+3)/2$.
In this case $\L$ is special if $d \leq 2m-2$
and $\L$ is non-special if $d=2m-1$ or $d=2m$.
\item[(f)] If $n = 2$ and $d \geq 2m+1$ then $\L$ is non-special and
$\l = dm+d-3m^2/2-m/2$.
\item[(g)] If $n \geq 2$ and $d \geq 2m$ then $\l=\l(d-m,d-2m,n-2,m)$.
\item[(h)] If $2 \leq m \leq d \leq 2m-1$ and $n \geq 3$
then $\L$ is empty (and therefore non-special),
so that $\l = -1$.
\end{itemize}
\end{lemma}

\begin{pf}
Statements (a) and (b), where we have either no base points
or simple base points, are trivial.
Statements (c), (d), (e), and (f),
where we have $3$ or fewer multiple points,
are handled easily by putting the $3$ points
at the coordinate points of the plane
and counting homogeneous monomials.
Statement (g) is obtained by making a quadratic Cremona transformation
at the point $p_0$ of multiplicity $d-m$ and two of the $n$ points
of multiplicity $m$; we note that the resulting linear system
is of the same form, namely that it is quasi-homogeneous with $m_0=d-m$.

Finally we turn to statement (h).
If $n \geq 3$, $d < 2m$ and $\L$ is nonempty
then the line $L_{ij}$ through any two of the $n$ points $p_i$ and $p_j$
must split off the linear system,
since $\L\cdot\L(1,0,2,1) = d-2m$.
If in fact $n \geq 4$, then the two lines $L_{12}$ and $L_{34}$ become
$(-1)$-curves on the blowup of the plane, which meet at one point;
hence the sum $L_{12}+L_{34}$ moves in a pencil,
and so cannot be part of the fixed part of the system $\L$.
This contradiction shows that $\L$ must be empty.

To finish we may therefore assume that $n=3$.
Again the three lines through the three points split off the system,
(in fact each splits off $2m-d$ times)
leaving the residual system $\L(4d-6m,d-m,3,2d-3m)$.
Therefore $\L$ is clearly empty if $2d<3m$.

If $2d \geq 3m$, then the $3$ lines
through $p_0$ and the $3$ points $p_i$
split off the residual system,
since the intersection number
$\L(4d-6m,d-m,3,2d-3m)\cdot\L(1,1,1,1)
= (4d-6m)-(d-m)-(2d-3m) = d-2m<0$.
Indeed, we see that these three lines each must split off $2m-d$ times,
leaving as the further residual system $\L(7d-12m,4d-7m,3,3d-5m)$
if all of these numbers are non-negative.
If $7d < 12m$ we see then that
the residual system, and hence $\L$, is empty, and we are done.
Note that if $7d \geq 12m$ then $3d \geq 5m$ (since $5/3<12/7$).

If $4d\leq 7m$ then the residual system is $\L(7d-12m,0,3,3d-5m)$;
again the $3$ lines through the $3$ points split off this system,
each $2m-d$ times, leaving the system $\L(10d-18m,0,3,5d-9m)$;
but if $4d\leq 7m$ then (since $7/4<9/5$) we have that $5d<9m$,
so that this further residual system is empty, and we are done.

If on the other hand $4d > 7m$ then the residual system is actually
$\L(7d-12m,4d-7m,3,3d-5m)$, and all these numbers are strictly positive.
Define the ratio $r = d/m$.
We have shown that if $r\leq 7/4$ then $\L(d,d-m,3,m)$ is empty,
while if $r > 7/4$ then $\l(d,d-m,3,m)=\l(7d-12m,4d-7m,3,3d-5m)$.
This residual system is of the same form ($m_0=d-m$) and has as its ratio
$s = (7d-12m)/(3d-5m) = (7r-12)/(3r-5)$.
Therefore if $s < 7/4$,
which happens for $r<13/7$, the system is again empty.

We claim that by iterating this procedure enough times we will be done,
i.e., for any $r < 2$ there is an iterate $s^{(n)}(r)$
which is less than $7/4$.
The function $s(r)$ maps the interval $r\in(7/4,2)$
onto the interval $s \in (1,2)$,
and moreover $s(r) < r$ for each such $r$.
Hence the iterates $s^{(n)}(r)$ form a decreasing sequence,
and if they never go below $7/4$,
they must converge to a fixed point of the function $s(r)$.
This is impossible, since the only fixed point is at $r=2$.
\end{pf}

The above Lemma allows the construction of an algorithm
to compute $\l(d,d-m,n,m)$.
If $n \geq 2$ and $d \geq 2m$ one uses (g) to reduce the numbers,
and hence one may assume that either $d < 2m$ or $n<2$.
Each of these cases is covered by the Lemma.
One can turn this algorithm into a formula,
and a criterion for speciality,
without too much difficulty.

\begin{proposition}
\label{m0=d-m}
Let $\L=\L(d,d-m,n,m)$ with $2 \leq m \leq d$.
Write $d = qm+\mu$ with $0\leq \mu\leq m-1$,
and $n = 2h+\epsilon$, with $\epsilon \in\{0,1\}$.
Then the system $\L$ is special if and only if
$q = h$, $\epsilon = 0$, and $\mu \leq m-2$.
More precisely:
\begin{itemize}
\item[(a)] If $q \geq h+1$ then $\L$ is nonempty and non-special.
In this case 
\[
\dim\L=d(m+1)-\binom{m}{2}-n\binom{m+1}{2}.
\]
\item[(b)] If $q=h$ and $\epsilon = 1$
the system $\L$ is empty and non-special.
\item[(c)] If $q=h$, $\epsilon = 0$, and $\mu = m-1$,
the system $\L$ is nonempty and non-special;
in this case
\[
\dim\L=(m-1)(m+2)/2.
\]
\item[(d)] If $q=h$, $\epsilon = 0$, and $\mu \leq m-2$,
the system $\L$ is special;
in this case
\[
\dim \L = \mu(\mu+3)/2.
\]
\item[(e)] If $q\leq h-1$
the system $\L$ is empty and non-special.
\end{itemize}
\end{proposition}

\begin{pf}
If $q \geq h+1$, we may apply quadratic Cremona transformations
$h$ times, arriving at the system $\L(d-hm,d-(h+1)m,\epsilon,m)$,
which is nonempty and non-special.
If $q \leq h$, we may apply quadratic Cremona transformations
$q-1$ times, arriving at the system $\L(\mu+m,\mu,2(h-q)+\epsilon+2,m)$.
If either $q < h$, or $\epsilon = 1$,
then we apply Lemma \ref{m0=d-m_algorithm}(h)
and conclude that $\L$ is empty, and therefore non-special.

We are left with the case $q=h$ and $\epsilon = 0$,
for which we have the system $\L(\mu+m,\mu,2,m)$;
we then apply Lemma \ref{m0=d-m_algorithm}(e)
to conclude the proof.
\end{pf}

This analysis of the $m_0=d-m$ case
applies immediately when $m_0 > d-m$ also.
Consider the system $\L(d,d-m+k,n,m)$ with $k \geq 1$.
We note that the $n$ lines through $p_0$ and $p_i$
split off, each $k$ times,
leaving as the residual system the system
$\L(d-kn,d-kn-m+k,n,m-k)$
(which is of the type discussed above).
The speciality of $\L$ is then deduced from this residual system:

\begin{corollary}
\label{m0=d-m+k}
Let $\L = \L(d,d-m+k,n,m)$ with $k \geq 1$,
and let
\[
\L' = \L(d-kn,d-kn-m+k,n,m-k).
\]
Then $\dim \L = \dim \L'$ and
$\L$ is non-special unless either
\begin{itemize}
\item[(a)]
$k \geq 2$ and $\L'$ is nonempty and non-special,
or
\item[(b)] $\L'$ is special.
\end{itemize}
\end{corollary}

Finally we turn to the case when $m_0 = d-m-1$.

\begin{proposition}
\label{m0=d-m-1}
Let $\L=\L(d,d-m-1,n,m)$ with $2 \leq m \leq d-1$.
Write $d = q(m-1)+\mu$ with $0\leq \mu \leq m-2$,
and $n = 2h+\epsilon$, with $\epsilon \in\{0,1\}$.
Then the system $\L$ is non-special of dimension 
$d(m+2)-(n+1)m(m+1)/2$
unless
\begin{itemize}
\item[(a)] $q=h+1$, $\mu=\epsilon=0$, and $(m-1)(m+2) \geq 4h$,
in which case
\[
\dim \L = (m-1)(m+2)/2 - 2h,
\]
or
\item[(b)] $q = h$, $\epsilon = 0$, and $4q \leq \mu(\mu+3)$,
in which case $\dim \L = \mu(\mu+3)/2 - 2q$.
\end{itemize}
\end{proposition}

\begin{pf}
First we note that performing a quadratic Cremona transformation
to the system $\L(d,d-m-1,n,m)$
gives the subsystem of $\L(d-m+1,d-2m,n-2,m)$
(which is of the same type, namely ``$m_0=d-m-1$'')
with two general simple base points.
Therefore by induction if we perform $k$ such transformations,
we obtain the subsystem of $\L(d-k(m-1),d-(k+1)(m-1)-2,n-2k,m)$
with $2k$ general simple base points,
if $d-(k+1)(m-1)-2 \geq 0$ and $n \geq 2k$.

If $q \geq h+2$, or if $q = h+1$ and $\mu \geq 2$,
then we may perform $h$ quadratic Cremona transformations
and arrive at the subsystem of $\L(d-h(m-1),d-(h+1)m-2,\epsilon,m)$
with $2h$ general simple base points.
This system is non-special.

We now analyze the case with $q=h+1$ and $\mu \leq 1$.
If $m=2$ the system is easily seen to be empty,
and therefore non-special;
hence we assume that $m \geq 3$.
We perform $h-1$ transformations,
leading to the subsystem of
$\L(d-(h-1)(m-1),d-h(m-1)-2,2+\epsilon,m)
=\L(2m-2+\mu,m-3+\mu,2+\epsilon,m)$
with $2h-2$ general base points.
We note that each line joining $2$ of the $2+\epsilon$ points
of multiplicity $m$ in this system splits off,
with multiplicity $2-\mu$;
hence if $\mu = 0$ and the system is not empty,
it is certainly special.

If $\mu = 0$ and $\epsilon = 1$,
when $m=3$ the original system is $\L(4,0,3,3)$,
which is empty; if $m \geq 4$
the residual system is $\L(2m-8,m-3,3,m-4)$.
Performing a Cremona transformation on this gives the obviously
empty system $\L(m-4,m-3)$.

If $\mu = 1$ and $\epsilon = 0$
the residual system is $\L(2m-2,m-2,2,m-1)$
which transforms to the non-special system $\L(m,0,2,1)$;
hence $\L$ is non-special in this case.

If $\mu=\epsilon = 1$,
then the residual system is $\L(2m-4,m-2,3,m-2)$
which transforms to $\L(m-2,0,1,m-2)$
which is again non-special.

If $\mu = \epsilon = 0$,
then the residual system is $\L(2m-4,m-3,2,m-2)$
which transforms to $\L(m-1,0,2,1)$ and is therefore non-special.
The original system is therefore nonempty if and only if
$(m-1)(m+2) \geq 4h$, and this is the only special case.

If $q \leq h$, we perform $q-1$ transformations,
arriving at the system $\L(\mu+m-1,\mu-2,2(h-q)+\epsilon+2,m)$,
with in addition $2q-2$ simple base points.
(If $\mu = 0$ the system is empty, and therefore non-special;
if $\mu = 1$, the system is also empty unless $q=h=1$ and $\epsilon=0$,
but this implies $d=m$ which is impossible.)

If $\mu \geq 2$, and $h > q$, then we argue as in the proof of
Lemma \ref{m0=d-m_algorithm}(h) and conclude that the system is empty.

We are left to analyze the case $\mu \geq 2$ and $h=q$.
If $\epsilon=1$, then the three lines
through the three points split off, each with multiplicity $m-\mu+1$;
therefore if $2m \geq 4\mu-3$, the system is empty
(the residual has negative degree).
Otherwise the residual system is $\L(4\mu-2m-4,\mu-2,3,2\mu-m-2)$
and which transforms to $\L(2\mu-m-2,\mu-2)$.
Since $\mu \leq m-2$, this system is empty.

Finally we take up the case where $\mu \geq 2$, $h=q$ and $\epsilon=0$,
in which case the line through the two remaining points
splits off $m+1-\mu$ times, leaving the residual system
$\L(2\mu-2,\mu-2,2,\mu-1)$ (with $2q-2$ simple base points).
We perform one more transformation,
giving the system of plane curves of degree $\mu$ with $2q$ simple base
points, leading to the last exception.
\end{pf}

\section{$(-1)$-Special Systems with $m \leq 3$}

Suppose $\L(d,m_0,n,m)$ is a $(-1)$-special system
with $m \leq 3$.
Then $\L$ must be of the form $\L = \M + NC$
for some $N = 2$ or $3$,
where $C\in\L(\delta,\mu_0,n,1)$
is either a quasi-homogeneous $(-1)$-curve
or a (compound) quasi-homogeneous $(-1)$-configuration,
and $\v(\M) \geq 0$ and $\M \cdot C = 0$.
This implies that $C$ is either:
\begin{center}
\begin{tabular}{cc}
$\L(2,0,5,1)$ & \\
$\L(e,e-1,2e,1)$ & with $e \geq 1$, or \\
$\L(e,e,e,1)$ & with $e \geq 1$.
\end{tabular}
\end{center}
The last one on the list is compound when $e \geq 2$;
all others are irreducible $(-1)$-curves.

These observations are sufficient to classify such systems.

\begin{lemma}
\label{obirreg23}
The quasi-homogeneous $(-1)$-special systems with $m \leq 3$
are the systems $\L(d,m_0,n,m)$ on the following list:
\begin{center}
\begin{tabular}{cccc}
$\L(4,0,5,2)$ & & $\v=-1$ & $\l = 0$ \\
$\L(2e,2e-2,2e,2)$, & $e \geq 1$ & $v = -1$ & $\l = 0$\\
$\L(d,d,e,2)$, & $d \geq 2e\geq 2$ & $\v = d-3e$ & $\l = d-2e$ \\
$\L(4,0,2,3)$ & & $\v=2$ & $\l = 3$ \\
$\L(6,0,5,3)$ & & $\v=-3$ & $\l=0$ \\
$\L(6,2,4,3)$ & & $\v=0$ & $\l=1$ \\
$\L(3e,3e-3,2e,3)$, & $e\geq 1$ & $\v = -3$ & $\l=0$ \\
$\L(3e+1,3e-2,2e,3)$, & $e\geq 1$ &  $\v=1$ & $\l=2$ \\
$\L(4e,4e-2,2e,3)$, & $e\geq 1$ & $\v = -1$ & $\l=0$ \\
$\L(d,d-1,e,3)$, &  $2d\geq 5e \geq 5$ & $\v=2d-6e$ & $\l = 2d-5e$ \\
$\L(d,d,e,3)$, & $d \geq 3e \geq 3$ & $\v = d-6e$ & $\l = d-3e$ .
\end{tabular}
\end{center}
\end{lemma}

\begin{pf}
We'll only present the analysis for $m = 3$;
the $m=2$ case is similar and easier, and we leave it to the reader.
Using the notation above, with $\L = \M + N C$,
in this case $N$ may be either $2$ or $3$.
We first discuss when $N=2$.

Suppose that $\L(2,0,5,1)$ splits twice off $\L(d,m_0,5,3)$.
Then $\M=\L(d-4,m_0,5,1)$, and $\M\cdot\L(2,0,5,1) = 2d-13$,
which can never be zero; hence this case cannot occur.

Suppose that $A=\L(e,e-1,2e,1)$ splits twice off $\L(d,m_0,2e,3)$.
Then $\M=\L(d-2e,m_0-2e+2,2e,1)$,
and
\begin{align*}
0 = \M\cdot A &= (d-2e)e-(m_0-2e+2)(e-1)-2e \\
&= de-m_0(e-1)-6e+2
&= e(d-m_0-6)+m_0+2
\end{align*}
so that certainly $m_0\geq d-5$.
Clearly for $\v(\M) \geq 0$ we must have $m_0\leq d-2$, but in this case
\begin{align*}
0 &= de -m_0(e-1) -6e +2 \\
&\geq de -(d-2)(e-1) -6e +2 \\
&= d-4e
\end{align*}
so that $d \leq 4e$.

We take up in turn the various possibilities for $m_0$.

If $m_0=d-5$, then $e=m_0+2$ for $\M\cdot A = 0$,
which gives that $m_0 = e-2$, a contradiction
since $m_0 \geq 2e-2$ and $e \geq 1$.

If $m_0=d-4$, then $\v(\M) = 3d-8e-1\geq 0$,
but then $0=\M\cdot A = m_0+2-2e=d-2-2e$, so that $d=2e+2$
and hence $3(2e+2)-8e-1\geq 0$, forcing $e\leq 2$.
When $e=1$ we have $\L(4,0,2,3)$, which has $\v = 2$ but $\l = 3$.
When $e=2$ we have $\L(6,2,4,3)$, which has $\v=0$ but $\l=1$.

If $m_0=d-3$, then $\v(\M) = 2d-6e$, so $d \geq 3e$;
then $0=\M\cdot A = m_0+2-3e = d-3e-1$, so that $d=3e+1$.
This gives the system $\L(3e+1,3e-2,2e,3)$,
which has $\v=1$ but after splitting off $A$ twice
leaves the system $\M=\L(e+1,e,2e,1)$,
which has $\l=2$.

If $m_0=d-2$, then $\v(\M) = d-4e$, so $d\geq 4e$,
and now $0 = \M\cdot A = m_0+2-4e$, so that $m_0=4e-2$,
forcing $d=4e$.  This gives the special system
$\L(4e,4e-2,2e,3)$, which has $\v = -1$
but consists of the fixed curve $2A$
and the residual system $\M = \L(2e,2e,2e,1)$
which is non-empty (it is a quasi-homogeneous $(-1)$-configuration,
consisting of the $2e$ lines through $p_0$ and $p_i$).

Finally suppose that the compound class $\L(e,e,e,1)$
splits twice off $\L(d,m_0,e,3)$.
Here $A = \L(1,1,1,1)$, and $\M=\L(d-2e,m_0-2e,e,1)$,
so that $0 = \M\cdot A = d-m_0-1$, forcing $m_0=d-1$.
Then $\v(\M) = 2d-5e$, leading to the special systems
$\L(d,d-1,e,3)$ with $2d\geq 5e$.

This completes the analysis for $m=3$ and $N=2$.
We now turn to the $m=N=3$ case.

Suppose first that $A=\L(2,0,5,1)$ splits three times off $\L(d,m_0,5,3)$.
Then $\M=\L(d-6,m_0)$, which has $\v\geq 0$ if $m_0\leq d-6$
and $\M\cdot A = 2d-12$, forcing $d = 6$, and $m_0=0$,
leading to the system $\L(6,0,5,3)$.

Suppose that $A=\L(e,e-1,2e,1)$ splits three times off $\L(d,m_0,2e,3)$.
Then $\M = \L(d-3e,m_0-3e+3)$ so that $d\geq m_0+3$ for $\v(\M)\geq 0$,
and $0=A\cdot\M = de-m_0(e-1)-6e+3 = e(d-m_0-6)+m_0+3$,
so that certainly $m_0\geq d-5$.

If $m_0=d-5$, then $e=m_0+3=d-2$;
but $e+2 = d\geq 3e$ forces $e=1$ and $m_0= -2$, a contradiction.

If $m_0=d-4$, then $m_0=2e-3$, and so $d=2e+1$.
Again since $d \geq 3e$ forces $e=1$ but then $m_0 = -1$,
a contradiction.

If $m_0 = d-3$, then $m_0 = 3e-3$, and so $d = 3e$,
so that $\M=0$ and we have the system $\L(3e,3e-3,2e,3)$.

Finally suppose that the compound class $\L(e,e,e,1)$
splits three times off $\L(d,m_0,e,3)$.
Here $A = \L(1,1,1,1)$, and $\M=\L(d-3e,m_0-3e)$,
so that $0 = \M\cdot A = d-m_0$, forcing $m_0=d$.
Then $\v(\M) = d-3e$, leading to the special systems
$\L(d,d,e,3)$ with $d\geq 3e$.

This completes the $m=N=3$ analysis.
\end{pf}

Since we have discussed in some detail the speciality
of systems $\L(d,m_0,n,m)$ with $m_0 \geq d-m-1$
in Section \ref{sec10},
we take the opportunity to make the following observation:

\begin{corollary}
Let $\L = \L(d,m_0,n,m)$ with $m_0 \geq d-m-1$
and $m \in \{2,3\}$.  Then $\L$ is special
if and only if $\L$ is $(-1)$-special.
\end{corollary}

\begin{pf}
Let us first discuss the $m=2$ case.
Suppose that $m_0 = d-3$.
Then Proposition \ref{m0=d-m-1} gives no special systems.
Suppose next that $m_0 = d-2$.
The Proposition \ref{m0=d-m} gives that the only such special system
is $\L(2e,2e-2,2e,2)$, which is $(-1)$-special.
If $m_0 = d-1$,
Corollary \ref{m0=d-m+k} gives no special systems.
Finally if $m_0 = d$, Corollary \ref{m0=d-m+k}
gives the assertion.

Now we turn to the $m=3$ case,
and first suppose that $m_0 = d-4$.
Then Proposition \ref{m0=d-m-1}(a) leads only to the systems
$\L(4,0,2,3)$ and $\L(6,2,4,3)$
which are $(-1)$-special.
Proposition \ref{m0=d-m-1}(b) gives no special systems with $m=3$.

Next suppose that $m_0=d-3$.
Then Proposition \ref{m0=d-m} gives that the only such special systems
are $\L(3e,3e-3,2e,3)$ and $\L(3e+1,3e-2,2e,3)$ for $e \geq 1$.

If $m_0 = d-2$,
Corollary \ref{m0=d-m+k} gives only the special systems
$\L(4e,4e-2,2e,3)$, for $e\geq 1$.

The other cases of $m_0 \geq d-1$ are trivial.
\end{pf}

\section{The Classification of Special Systems with $m \leq 3$}

\begin{theorem}
\label{thmreg2}
A system ${\L}(d,m_0,n,m)$ with $m \leq 3$ is special
if and only if it is a $(-1)$-special system,
i.e., it is one of the systems listed in Lemma \ref{obirreg23}.
\end{theorem}

\begin{pf}
We will outline the proof in the case $m=3$;
the $m=2$ case is analogous in every way.

We will assume that $m_0 \neq 1,3$.
We may also assume $\L$ is not empty,
otherwise it is certainly non-special.
We will prove the theorem by induction on $n$;
the assertion is true for $n\leq 2$.
So we may assume $n\geq 3$.
The theorem is easily seen to be true for $d\leq 5$.
So we will assume $d\geq 6$.

Furthermore the cases $d<m_0$ are trivial
and the cases $d\leq m_0+4$ are taken care of by the results
of Section \ref{sec10}.
Hence we may assume that $d \geq m_0+5$, or that $m_0 \leq d-5$.

The general approach is to assume $\L$ is not $(-1)$-special
and prove that it is non-special.
We proceed by induction and we
assume the theorem holds for lower values of $n$.

We start with the case ${\v}\leq -1$.
We perform a $(2,b)$-degeneration and we require
that both kernel linear systems $\LPhat$, $\LFhat$ are empty.
The requirement that
$\LFhat={\L}(d,d-1,b,3)$ is empty translates into the inequality $5b>2d$.
As for requirement
that $\LPhat={\L}(d-3,m_0,n-b,3)$ is empty,
we can use induction and impose that its virtual dimension
${\vPhat}$ is negative, unless
$\LPhat$ is a $(-1)$-special system. This might
only happen (when $m_0 \leq d-5$) if:
(i) $m_0=d-5$, $d-3=4e$, $n-b=2e$;
(ii) $m_0=d-6$, $d-3=3e,3e+1$, $n-b=2e$;
(iii) $d=9, m_0=0, n-b=5$;
(iv) $d=9, m_0=2, n-b=4$;
(v) $d=7, m_0=0, n-b=2$.

Since ${\v}(d,m_0,n,3)\leq -1$,
the condition ${\v}(d-3,m_0,n-b,3)<0$ is implied by
${\v}(d-3,m_0,n-b,3)\leq {\v}(d,m_0,n,3)$,
which reads $2b\leq d$.
In conclusion we need to choose a $b$ such that
$\frac{2d}{5} < b \leq \frac{d}{2}$.
For later purposes we will need  $b < \frac{d}{2}$
if $d$ is divisible by $4$.
Since $\frac{d}{2} - \frac{2d}{5} = \frac{d}{10}$
we see that we can choose a suitable $b$ as soon as $d\geq 11$,
but in fact one sees that for $d$ between $6$ and $10$
and unequal to $8$,
the maximum $b$ such that $2b \leq d$ works.
On the other hand one proves directly the theorem for $d=8$,
so we can dispense with this case.

With this choice,
only the $(-1)$-special systems
with negative virtual dimension could occur, i.e.
(i),
(ii) with $d=3e+3$,
(iii).
We also remark that, with this choice of $b$, one has $b<n$.

Suppose we are in case (i).
Then $d=4e+3$ and if $d \geq 19$,
then we may choose $b$ in at
least two ways, and avoid the $(-1)$-special systems.
We are left with the cases $d=7,11,15$.
Then $b=3,5,7$ respectively, and,
if $\LPhat$ is a $(-1)$-special system,
then $n=5,9,13$ respectively.
But in each one of these cases the virtual dimension of $\L$ is $2$,
contrary to the hypothesis.

Suppose we are in  case (ii), with $d=3e+3$.
Again we can choose $b$ in at
least two ways, and avoid the $(-1)$-special systems,
as soon as $d\geq 18$.
The remaining cases $d=6,9,12,15$,
in which the values of $b$ are $b=3,4,5,7$ respectively
and therefore $n=5,8,11,15$,
also have non-negative virtual dimension for $d \geq 9$;
the case $d=6$ is a $(-1)$-special system.

In case (iii),
$b=4$ and therefore $n=9$,
in which case $\L$ has virtual dimension $0$, a contradiction.

We have therefore arranged to choose a $b$
in such a way that both kernel linear systems $\LFhat$ and $\LPhat$
are empty, unless $\L$ is $(-1)$-special.

Now we claim that, with the choices we made,
${\Lo}$ is empty, which implies that also $\L$ is empty,
hence non-special.
The assertion is clear if either one of the two systems
$\LP$, $\LF$ is empty.

So we may assume that both systems $\LP$ and $\LF$ are not empty.
The system $\LF$ is ${\L}(d,d-2,b,3)$ and ${\vF}=3d-1-6b$.
The $b$ lines through $p_0$ split off this system,
and the residual system is $\L(d-b,d-b-2,b,2)$,
which by Theorem \ref{thmreg2} and Lemma \ref{obirreg23}
could only be special if $d=4e$ and $b=2e$.
Because of our choice of $b$,
this does not occur; hence $\LF$ is non-special
and therefore its dimension is ${\lF}={\vF}=3d-1-6b > -1$.

Let us examine the possibility that $\LP=\L(d-2,m_0,n-b,3)$
is a $(-1)$-special system.
For $d \geq 9$, or $d=6$ or $7$, this can only be the case if:
(i) $m_0=d-5$, $d=3e+2,3e+3$, $n=b+2e$;
(ii) $d=6, m_0=0, n=b+2$.
In case (ii) $\L$ is the $(-1)$-special system $\L(6,0,5,3)$.
In cases (i) with $d=3e+2$, then ${\rP} = \lP =0$;
since $b>0$, the linear series $\RF$ contains $b$ general fixed points,
and hence $\RF \cap \RP$ is empty.
Therefore by Lemma \ref{dimLo2}
${\Lo} = -1$.
Similarly in case (i), if $d=3e+3$,
then ${\lP}=2$ but $b>2$, which again forces $\RF \cap \RP$ to be empty.

Hence we may now assume that also $\LP$ is non-special.
Now we conclude that $\L$ is empty by Corollary \ref{cordimLk}(a).

Now we consider the case ${\v}\geq 0$.
Again we may assume $m_0 \leq d-5$.
We observe that we may assume $n>\frac{5}{6}d-1$,
otherwise ${\v}(d,d-4,n,3)\geq 1$
and therefore ${\L}(d,d-4,n,3)$ is non-special,
hence also ${\L}(d,m_0,n,3)$ is non-special.

We propose to perform a $(3,h)$-degeneration,
where we write $d = 2h-\epsilon$ with $\epsilon \in \{0,1\}$,
With this choice $\LFhat=\L(d,d-2,h,3)$ is empty
if $d$ is not divisible by $4$,
and has dimension $0$ if $4$ divides $d$,
since $\l(d,d-2,h,3)=\l(d-h,d-h-2,h,2)$
and $\v(d-h,d-h-2,h,2) = 3d-6h-1 = -3\epsilon-1 < 0$,
and is non-special unless it is $(-1)$-special,
which can only happen if $h=2e$ and $d=2h$.

Now $\LF = \L(d,d-3,h,3)$ which is never $(-1)$-special;
hence $\lF = \v(d,d-3,h,3) = 4d-3-6h = d-3(\epsilon+1)$.
Therefore $\rF = d-3(\epsilon+1) - \eta$
where $\eta = 1$ if $4|d$ and is zero otherwise.

Now note that
$\vP-\v = \v(d-3,m_0,n-h,3)-\v(d,m_0,n,3) = 6h-3d = 3\epsilon \geq 0$,
so that $\LP$ is nonempty.
Suppose first that $\LP$ is non-special,
and $\LPhat$ is empty.
Then we have $\rP = \vP = \v+3\epsilon$,
and so $\rF+\rP = d-3 - \eta + \v \geq d-4$.
Hence we apply Proposition \ref{dimLo3}(b)
and finish.

Suppose next that $\LP$ is non-special,
and $\LPhat$ is nonempty and non-special.
Then $\rP=\vP-\vPhat-1 = d-3$ so that Proposition \ref{dimLo3}(b)
applies and we finish.

Next we suppose that $\LP$ is non-special
but $\LPhat$ is nonempty and special.
In all the possibilities for special $\LPhat$ with $d \geq 9$,
except one,
we have $\lPhat \leq 2$;
hence $\rP \geq d-6$, which forces $\rP+\rF \geq d-4$ easily,
and we again finish using Proposition \ref{dimLo3}(b).
The one exception to this is when $m_0 = d-5$
and $2d-8 \geq 5(n-h)$.
However in this case $\lPhat = 2(d-4)-5(n-h)$,
and since $\v \geq 0$, we see easily that $\rP +\rF \geq d-4$ again.

Finally we must discuss the case that
${\LP}$ is $(-1)$-special.
Since we have already shown that $\vP \geq 0$,
this can only be the case (given that $d\geq 9$) when:
(i) $d=9$, $m_0=2$, $h=5$, $n = 9$;
(ii) $d=3e+4$, $m_0=3e-2$, $n=h+2e$.

In the first case we see that $\v = \v(9,2,9,3) = -3$,
a contradiction.

As to the second case, we suggest instead a $(3,h+1)$-degeneration.
In this case $\LFhat$,
which had dimension at most $0$ before,
is now empty.
Moreover $\LF = \L(d,d-3,h+1,3)=\L(3e+4,3e+1,(3e+6+\epsilon)/2,3)$
is again never $(-1)$-special.
Hence $\rF = \lF = \vF = 3e-5-3\epsilon$.

Now $\LPhat = \L(3e,3e-2,2e-1,3)$ which is never $(-1)$-special;
since $\vPhat = 5-3e < 0$, we see that $\LPhat$ is empty.

The system $\LP = \L(3e+1,3e-2,2e-1,3)$ is also never $(-1)$-special;
moreover $\vP = 7$, so that $\LP$ is not empty and also $\rP = 7$.

Therefore $\rF+\rP = 3e+2-3\epsilon = d-2-3\epsilon$;
hence if $\epsilon = 0$, we conclude the proof using
Proposition \ref{dimLo3}(b).

However if $\epsilon = 1$, then $e$ is odd, and $h = (3e+5)/2$,
so that $n = (7e+5)/2$;
computing the virtual dimension of the system $\L$
we find $\v = -2$, a contradiction.

This completes the proof in the $m=3$ case.

As noted above, we leave the details of the $m=2$ case to the reader.
The outline of the proof is to make a $(1,b)$ degeneration,
and again use induction on $n$.
We may assume that $d \geq m_0+4$.

When ${\v}\leq -1$,
we need to prove that if $\L$ is not $(-1)$-special,
then it is empty.
We perform a $(1,b)$-degeneration,
where $b$ is the minimum integer with $2b>d$;
the relevant linear systems are
$\LP = \L(d-1,m_0,n-b,2)$,
$\LF = \L(d,d-1,b,2)$,
$\LPhat = \L(d-2,m_0,n-b,2)$, and
$\LFhat = \L(d,d,b,2)$.

The reader can check that with this choice of $b$,
both $\LFhat$ and $\LPhat$ are empty.
Hence if we can show that $\lP+\lF \leq d-2$,
then we can apply Corollary \ref{cordimLk}(a)
and conclude that $\L$ will be empty.
This is automatic if either $\LP$ or $\LF$ is empty,
so one may assume that both systems are not empty.
If so, one checks as in the $m=3$ case
that we do have $\lP+\lF \leq d-2$ as needed.

In case ${\v}\geq 0$,
still assuming that $d \geq 4$ and $m_0 \leq d-3$.
We notice we can also assume $n\geq d$;
otherwise, as we proved already,
${\L}(d,d-2,n,2)$ is non-special and not empty,
and therefore also ${\L}(d,m_0,n,2)$ is non-special
for all $m_0\leq d-3$ by Lemma (0.1)\ref{basic1}(c).

Again we perfom a $(1,b)$-degeneration,
and we take $b$ to be the maximum such that $2b \leq d+1$.
Then $\lFhat = \l(d,d,b,2) = d-2b$,
and $\LF$ has dimension $\lF = 2d-3b$.
Therefore the system $\LF$ is non-special,
and $\rF = \lF-\lFhat-1 = (2d-3b)-(d-2b)-1 = d-b-1$.

One checks that $\LP$ is non-special and non-empty.
Since we already know that $\LF$ is non-special and non-empty,
we will be done if we show that $\rP+\rF \geq d-2$,
by applying Proposition \ref{dimLo3}(b),
and noting that in this case
$\lo=\lP+\lF-d+1=\vP+\vF-d+1=\v$.

Since we have seen above that $\rF=d-b-1$,
we need only to show that $\rP \geq b-1$ to finish the proof.
This we leave to the reader.
\end{pf}

\end{document}